\documentclass[12pt]{amsart}
\usepackage[utf8]{inputenc}
\usepackage{amsfonts, amssymb, euscript,graphics}
\usepackage{amsmath}
\usepackage{amscd}
\usepackage{graphicx}
\usepackage{eepic}
\usepackage{comment}
\usepackage{psfrag}
\usepackage{wrapfig}

\newcommand{\R}{{\mathbb R}}

\newcommand{\Z}{{\mathbb Z}}

\renewcommand{\em}[1]{{\par\noindent\bf #1}}

\newtheorem{theorem}{Theorem}[section]
\newtheorem{lemma}[theorem]{Lemma}

\newfont{\cmbb}{cmex10 at 7pt}

\newcommand{\ign}[1]{{}}
\def \.{\mskip 1mu}
\def \?{\mskip -1mu}

\renewcommand*\proof[1]{{\noindent\it Proof#1.\/\enspace\ignorespaces}}
\newcommand*\pr{\proof{}}

\def\sq{\raisebox{-.1ex}{\hbox{$\square$}}}

\def\qed{{\unskip\nobreak\hfill\hskip.5em\nobreak\hfil\sq
\parfillskip=0pt\medbreak}}

\makeatletter

\renewcommand{\abovecaptionskip}{.8ex}

\newcommand*{\capdelimiter}{.}
\renewcommand{\@makecaption}[2]{%
\vspace{\abovecaptionskip}%
\sbox{\@tempboxa}{#1\capdelimiter#2}
\ifdim \wd\@tempboxa >\hsize
   #1\capdelimiter#2 \par
\else
 \global\@minipagefalse \hbox to \hsize {\hfil #1\capdelimiter#2 \hfil}%
\fi
\vspace{\belowcaptionskip}}

\newcommand{\subs}[1]
   {\refstepcounter{subsection}
   \medskip\noindent
   {\it\arabic{section}.\arabic{subsection}.\hspace{.25em}\ignorespaces#1.}
          \ignorespaces}

   \@addtoreset{equation}{section}

\makeatother

\title{Chekanov-type theorem for spherized cotangent bundles}

\author{Petr E. Pushkar\\\today}
\address[Petr E. Pushkar]{National Research University Higher School of Economics (HSE)}
\email{petya.pushkar@gmail.com}
\thanks{The work  were partially supported by the grant RFBR 15-01-05990.}

\date{}

\begin{document}
\begin{abstract}
We prove a Chekanov-type theorem for the spherization of the cotangent bundle $ST^*B$ of a closed manifold $B$. It claims that for Legendrian submanifolds in $ST^*B$  the property  ``to be given by a generating family quadratic at infinity'' persists under Legendrian isotopies. 
\end{abstract}
\maketitle

\section*{Introduction} The notion of a generating family  is well known in symplectic and contact topology. The key role in many problems plays the theorem,  claiming that the property of Lagrangian (correspondingly Legendrian) submanifolds in $T^*B$ (corr. $J^1(B,\R)$) to be given by generating family is an invariant under Hamiltonian (corr. contact) isotopy (see, for example, \cite{Cha, Ch, EG, Fe, L-S, Vi}).    

In this paper we prove a similar result for Legendrian submanifolds in spherizations of cotangent bundles:
\begin{theorem}
\label{chst^*}
%
Let $B$ be a closed manifold, $E \to B$ be a smooth compact fibration. Let $\{L_t\}_{t \in [0,1]}$ be a legendrian isotopy of a compact Legendrian manifold $L_0 \subset ST^*B$. 
Suppose $L_0$ is given by a generating family $F\colon E \to \R$.
Then there exists $N\in \Z_+$, such that $L_t$ is given  by a generating family $G_t\colon E\times \R^N \to \R$ of the form:
$$
G_t(e,q)=F(e)+Q(q)+f_t(e,q)
$$
for a nondegenerate quadratic form $Q$ on $\R^N$ and compactly supported function $f_t$ such that $f_0=0$.
\end{theorem}

A very closed statement to the theorem \ref{chst^*} is the  Theorem 4.1.1 contained in \cite{EG} (see also \cite{Fe} for another close result) and the main aim of this paper is to fill a gap in the proof (see \cite{P1}).

\section{Generating families for Legendrian submanifolds in $ST^*$}

\subs{Contact structure on $ST^*$} We recall standard notions from contact geometry \cite{AG}. Let $B$ be a smooth manifold. 
Denote by $0_B$ the zero section of the cotangent bundle $T^*B$. 
The spherisation $ST^*B$ of the cotangent bundle $T^*B$ is a quotient space under the natural free action of a multiplicative group of positive real numbers $\R_+$ on $T^*B\setminus 0_B$, a positive number $a$ transform a pair $(q,p)$ $(q\in B, p \in T^*_qB)$ into the pair $(q, ap)$. 
The space $ST^*B$ is a smooth manifold of the dimension $2\dim B -1$ and it carries a natural cooriented contact structure $\xi$ defined as follows. 
Consider the  Liouville $1$-form $\lambda=pdq$ on $T^*B$, it defines a cooriented (by $\lambda$ itself) hyperplane field $\lambda=0$ on $T^*B\setminus 0_B$. 
This cooriented hyperplane field is invariant with respect to the action of $\R_+$ and it is tangent to orbits of the action. 
Hence the projection of that hyperplane field to $ST^*B$ is a cooriented hyperplane field on $ST^*B$ which turns out to be a contact structure.
 We remark here, that there is no natural choice of a contact form on $ST^*B$.

\subs{Critical points and critical values of Legendrian manifolds in $J^1(B)$} Consider the space $J^1(B)=J^1(B,\R)=T^*B\times\R$ of 1-jets of functions on a closed manifold $B$. 
We say that a point of a Legendrian manifold $\Lambda\subset J^1(B)$ is a critical point of $\Lambda$ if it projects to a zero section $0_B$ under the natural projection $J^1(B) \to T^*B$. 
We say that a number $a$ is a critical value of $\Lambda\subset J^1(B)$ if $a$ equals to the value of $u$-coordinate of a critical point of $\Lambda$ . 

\subs{Legendrian manifolds} The following notion generalise a notion of regular level set of a function on a manifold. Let us fix a number $c \in \R$.
Suppose that Legendrian manifold $\Lambda\subset J^1(B)$ is transverse to  $T^*B\times \{c\} \subset T^*B \times \R =J^1(B)$. 
Note that this implies that $c$ is not a critical value of $\Lambda$.
Consider a manifold $L^c = \Lambda \cap (T^*B\times \{c\})$. 
The intersection of $L^c$ with $0_B\times \{c\}$ is empty and the restriction of the natural projection $(T^*B\setminus 0_B)\times \{c\} \to ST^*B$ to $L^c$ is a legendrian immersion.
We denote the image of $L^c$ by $\Lambda^c$ and we say in this situation that $\Lambda^c$ is a $c$-reduction of $\Lambda$. 
Note that $ST^*B$ is a contact reduction (in the sence of \cite{EG,P1}) of $J^1(B)\setminus (0_B\times\R)$ along any hypersurface given by equation $u=c$ in $J^1(B)\setminus (0_B\times\R)$ and $\Lambda^c$ is a reduction of $\Lambda$.

\subs{Legendrian manifolds and generating families} 
Let us recall firstly the definition of generating family in the space of $1$-jets of functions.
Let $B$ be a manifold, consider the space $J^1(B)=T^*B\times \R$ of one jets of functions on $B$. 
The space $J^1(B)$ is a contact manifold with the canonical contact structure given by the form $du-\lambda$, where $\lambda$ is (a lift of) the Liouville form on $T^*B$, $u$ is the coordinate on the factor $\R$. 
A smooth bundle $E\to B$ and a generic function $F\colon E\to\R$ generates an (immersed) Legendrian submanifold $\Lambda_F \subset J^1(B)$ as follows. 
Consider a fiber of the bundle $E\to W$, and a critical point of the restriction of the function $F$ to this fiber. 
Denote by $C_F$ the set consisting of all such points.
For a sufficiently generic function $F$ the set $C_F$  is a smooth submanifold of the total space $E$.
 The genericity condition is that the equation $d_wF=0$, where $w$ is a local coordinate on a fiber (for a local trivialization of the bundle) of $E\to B$, satisfies the condition of the implicit function theorem. 
At any point $z$ of $C_F$ the differential $d_BF(z)$ of the function $F$ along the base $B$ is well defined. 
The rule $z\mapsto (z,d_BF(z),F(z))$ defines an immersion $l_F\colon C_F \to J^1(B)$ and its image is a Legendrian manifold $\Lambda_F$ under definition.
Function $F$ is called a generating family for $\Lambda_F$.
For a generating family $F$ in a local trivialization $B\times W$ of $E$ the manifold $\Lambda_F$ is given by the formula:  
$$
\Lambda_F = \{(q,p,u)| 
\exists w_0 F_w(q,w_0)=0, p = F_q(q,w_0), u=F(q,w_0)\},
$$
where $q,p$ are canonical coordinates on $T^*B$.

Now we define a Legendrian (immersed) submanifold $\Lambda^c_F$ in the space $ST^*B$ starting from a smooth bundle $\pi \colon E\to B$ and a generic function $F \colon E \to \R$. 
For a point $b\in B$ we denote by $C^c_F(b)$ the $c$-level of the restriction of the function $F$ to $C_F$. 
Genericity conditions are the following -- $C_F$ is a manifold in a neighborhood of $C^c_F$ and $c$ is a regular value of $F$ restricted to $C_F$. 
Consider a map $C^c_F \to ST^*B$: $z\mapsto (z,[d_BF(z)])$. 
This map is well defined since for $z \in C^c_F$ $d_BF(z)\ne 0$.
Moreover, this map is a legendrian immersion and we denote its image by $\Lambda^c_F$. 
We will be interested in embedded Legendrian manifolds only.

We remark here that if $F$ is a generating family for a manifold $\Lambda_F \subset J^1(B)$ then $F$ is a generating family for a manifold $L_F\subset ST^*B$ 
if and only if $\Lambda_F$ is transversal to the hypersurface $\{u=0\}$ in $J^1(B)$.

\subs{Stabilization} Similarly to the $J^1(B)$-case one can stabilize generating families for Legendrian manifolds in $ST^*B$. If $F\colon E \to \R$ is a generating family for a Legendrian $L=L_F\subset ST^*B$, then a function $G\colon E\times \R^K$ of the form $$G(e,w)=G(e)+Q(w), e\in E, w\in \R^K$$ where $Q$ is a non-degenerate quadratic form on $\R^k$  is also a generating family for $L$. 

\section{Legendrian isotopy lifting}

{\it Proof of Theorem \ref{chst^*}.} We follow the strategy of the proof of Theorem 4.1.1. in \cite{EG}.
The (generalized) Chekanov theorem \cite{Ch,P} for $J^1(B)$ has almost the same statement as Theorem \ref{chst^*} -- $ST^*B$ is replaced by $J^1(B)$. 
We can perturb the generating family $F$ to a generating family $\widetilde{F}$, such that $\widetilde{F}$ is a generating family for $L_0$ and at the same moment is a generating family for a closed submanifold in $J^1(B)$.  Now Theorem~\ref{chst^*} is a direct corollary of the corresponding $J^1(B)$-result \cite{P} in view of the following legendrian isotopy lifting lemma.
\qed 

\subs{Legendrian isotopy lifting lemma} The following lemma claims contact isotopy lifting property formulated in \cite{EG} for more general situations. Unfortunately the general statement is not true \cite{P1} and, correspondingly, our proof use the specificity of the $ST^*$-situation.

\begin{lemma} 
Let $B$ be a closed manifold and $c\in \R$. 
Consider a compact Legendrian manifold $\Lambda \subset J^1(B)$ such that its $c$-reduction $\Lambda^c$ 
is well defined and $\Lambda^c$ is an embedded manifold. 
Let $L_{t,t \in [0,1]}$ be a legendrian isotopy of $\Lambda^c=L_0$. 
Then there exists a legendrian  isotopy  $\Lambda_{t,t \in [0,1]}, \Lambda_0=\Lambda$ 
such that for any $t\in [0,1]$ its $c$-reduction is defined and $\Lambda^c_t=L_t$. 
\end{lemma}

\pr It is sufficient to prove the statement of the lemma for $c=0$. 
Consider the legendrian isotopy ${L_t}$. 
By isotopy extension theorem there exists a contact flow $\varphi_{t\in [0,1]}$, such that
$\varphi_t(L_0)=L_t$  for any $t \in [0,1]$.
Any contact isotopy of $ST^*B$ lifts to a (homogeneous) Hamiltonian flow on $T^*B\setminus 0_B$.
More precisely -- consider a Hamiltonian $H_t \colon T^*B\setminus 0_B \to \R$ such that $H_t(ap,q)=aH_t(p,q)$ for any positive number $a$ (we will say that such a Hamiltonian is homogeneous).
Then the flow of such a Hamiltonian function is well defined for all values of $t$ and projects to a contact flow on $ST^*B$. 
Moreover, any contact flow on $ST^*B$ could be given as a projection of a unique Hamiltonian flow above.

We take a homogeneous Hamiltonian $H_t$ corresponding to the flow~$\varphi_t$.
Consider a function $K_t(p,q,u)=H_t(p,q)$ on $(T^*B\setminus 0_B)\times\R \subset J^1(B)$ as a contact Hamiltonian (see \cite{AG}) with respect to the contact form $du-\lambda$.
Any set $(T^*B\setminus 0_B)\times \{c\}$ is invariant under the flow generated by $K_t$ and coincides on it with the flow of $H_t$ under the forgetful identifications $T^*B\times \{c\} = T^*B$. 
Indeed, it follows from the explicit formula for the corresponding contact vector field: $\dot{u}=K-pK_p, \dot{p}=K_q-pK_u,\dot{q}=-K_p$ (see \cite{AG}). 
The $u$-component of this contact vector field equals to zero since $K$ is homogeneous. 
Hence the flow $\psi_t$ generated by $K$ satisfy $\psi_t(\Lambda \cap (T^*B\setminus 0_B)=L_t$.
In general it is impossible to extend $\psi_t$ to a flow on the whole space $J^1(B)$ so we will change the function $K_t$.
Let us fix an arbitrary smooth function $\widetilde{H}_t\colon T^*B \to \R$ coinciding with $H_t$ in a neighbourhood of infinity.
Denote by $P_t$ the function $P_t(p,q,u)=\widetilde{H}_t(p,q)$ and by $P_t^C$ ($C\in \R_+$) the function $P_t^C(p,q,u)=\frac{1}{C}P_t(Cp,q,u)$.
We claim that for sufficiently big $C$ the legendrian isotopy of $\Lambda$ generated by the contact flow $\Psi_t^C$ of $P_t^C$ 
satisfies the claim of lemma. 

Let us fix a number $a$ such that the absolute value of any critical value of $\Lambda$ is bigger then $2a$. 
Denote by $X \subset \Lambda$ the subset formed by all points such that the absolute value of $u$-coordinate is at most $a$, by $Y$ we denote the closure of its complement $\Lambda\setminus X$.
The set $X$ is a compact set and contained in $(T^*B\setminus 0_B)\times \R$.
Take a neighborhood $U\subset T^*B$ of the zero section,  the support of $P_t^C-K_t$ is containes in  $U\times\R$ for sufficiently big $C$.
Hence, for sufficiently big $C$, $\Psi_t^C(X)=\psi_t(X)$ for all $t\in [0,1]$.
It remains to show that for sufficiently big $C$ the $u$-coordinate of any point in $\Psi_t^C(Y)$ could not be zero and hence zero reduction of $\Psi_t^C(\Lambda)$ is $L_t$.
The coordinate $u$ changes under the action of a contact Hamiltonian $P^C$ according to the low: $\dot{u}=P_t^C - p\frac{\partial P_t^C}{\partial p}$.
So it is sufficient to show that the speed of the $u$-coordinate uniformly tends to zero as $C$ tends to infinity. 
The following general consideration finishes the proof.

Consider a smooth vector bundle $V$ over a closed manifold $M$. We denote by $M(c)$ fiberwise multiplication by $c$.
We say that a smooth function on $V$ is positively homogeneous degree $1$ at infinity if
it is coincides with a continuous positively homogeneous (i.e. $1/\alpha (M(\alpha))^*$-invariant for any positive $\alpha$) degree $1$  function up to a sum with a compactly supported continuous function.
Let $v$ be a  vector field tangent to fibers of V and coinciding with Euler vector field on each fiber of~$V$. 
Consider an operator $D$ sending a function $g$ on $V$ to $g-L_vg$. 
For a positively homogeneous function $f$ the function $Df$ is a compactly supported function.
Denote by $f^C$ the function $1/C (M(C))^*f$, i.e. for any $x\in V$ $f^C(x)=1/Cf(M(C)x)$.

\begin{lemma} For any smooth positively homogeneous of degree $1$ at infinity function $f$, the $C^0$-norm of $D(1/C (M(C))^*f)$ tends to zero while $C\to +\infty$.
\end{lemma}

\pr Indeed, $D(\frac{1}{C}(M(C))^*f)=\frac{1}{C}(M(C))^*D(f)$. Hence the $C^0$-norm of $D(1/C (M(C))^*f)$ equals to the  $C^0$-norm of $f$ divided by $C$.
\qed


\begin{thebibliography}{99}

\bibitem[AG]{AG} V.I. Arnold, A. Givental, {Symplectic geometry},  Encyclopedia of mathematical sciences vol 4., Springer 1997.

\bibitem[Cha]{Cha} M. Chaperon, On generating families, in: The Floer Memorial Volume, Helmut Hofer ed., Progr. Math. 133, Birkhauser, 1995.

\bibitem[Ch]{Ch} Chekanov, Yu. V. Critical points of quasifunctions,
and generating families of
Legendrian manifolds.
Funct. Anal. Appl. { 30}
(1996), no. 2, 118--128.

\bibitem[F]{Fe} E. Ferrand, On a theorem of Chekanov, from: “Symplectic singularities and geometry
of gauge fields (Warsaw, 1995)”, Banach Center Publ. 39, Polish Acad. Sci., Warsaw
(1997) 39–48 

\bibitem[EG]{EG} Ya. Eliashberg, M. Gromov, Lagrangian intersection theory: finite-dimensional approach, Geometry
of differential equations,  27–118, Amer. Math. Soc. Transl. Ser. 2, 186, AMS, Providence, RI, 1998.

\bibitem[L-S]{L-S}Laudenbach, F., and Sikorav, J.-C.. "Persistance d'intersection avec la section nulle au cours d'une isotopie hamiltonienne dans un fibre cotangent.." Inventiones mathematicae 82 (1985): 349-358. 

\bibitem[P]{P} Pushkar', P.~E. A generalization of Chekanov's
theorem. Diameters of immersed manifolds and wave fronts. Proc. Steklov Inst. Math. 1998, no. 2 (221), 279--295.

\bibitem[P1]{P1} Pushkar', P.~E.  On Hamiltonian and contact isotopy liftings,
Preprint {\tt 	arXiv:1602.07948} 

\bibitem[Vi]{Vi} C. Viterbo, Symplectic topology as the geometry of
generating functions, {\em Math. Ann. } {\bf 292} (1992), 685-710.

\end{thebibliography}
\end{document}